\documentclass[11pt]{amsart}
\usepackage{epsfig}
\usepackage{graphics}

\newtheorem{theorem}{Theorem}
\newtheorem{lemma}[theorem]{Lemma}

\newtheorem{corollary}[theorem]{Corollary}

\theoremstyle{definition}

\newtheorem{remark}[theorem]{Remark}

\theoremstyle{remark}

 \def\Z{{\mathbb{Z}}}

\def\mcg{{\rm Mod}}

\begin{document}

\newenvironment{prooff}{\medskip \par \noindent {\it Proof}\
}{\hfill
$\square$ \medskip \par}
    \def\sqr#1#2{{\vcenter{\hrule height.#2pt
        \hbox{\vrule width.#2pt height#1pt \kern#1pt
            \vrule width.#2pt}\hrule height.#2pt}}}
    \def\square{\mathchoice\sqr67\sqr67\sqr{2.1}6\sqr{1.5}6}
\def\pf#1{\medskip \par \noindent {\it #1.}\ }
\def\endpf{\hfill $\square$ \medskip \par}
\def\demo#1{\medskip \par \noindent {\it #1.}\ }
\def\enddemo{\medskip \par}
\def\qed{~\hfill$\square$}

 \title[Homomorphisms from mapping class groups]
 {Homomorphisms from mapping class groups}

 \author{William Harvey}
\author{Mustafa Korkmaz}

 \address{Department of Mathematics, King's College, London, WC2R
2LS.}
\email{bill.harvey@kcl.ac.uk}

 \address{Department of Mathematics, Middle East Technical
University, 06531 Ankara, Turkey}
\email{korkmaz@arf.math.metu.edu.tr}

 \date{\today}

\begin{abstract}This paper concerns rigidity of the mapping class groups.
 We show that any homomorphism
$\varphi:\mcg_g\to\mcg_h$ between mapping class groups of
 closed orientable surfaces with distinct genera $g>h$
 is trivial if $g\geq 3$ and has finite image for all $g\geq 1$. Some implications
 are drawn for more general homomorphs of these groups.
\end{abstract}

 \maketitle
  \setcounter{secnumdepth}{1}
 \setcounter{section}{0}

\section{Introduction}
The mapping class groups resemble lattices in higher rank Lie
groups in various ways. Along these lines, in parallel with the
classic rigidity theorems of Mostow and Margulis, there are results
by various authors restricting the existence of nontrivial morphisms
between lattices and mapping class groups; for a discussion of
this the reader may consult the extensive review article by N.V.
Ivanov \cite{iv}.

There are natural inclusion morphisms between mapping class groups
of finite-type surfaces stemming from inclusions of bordered
surfaces, and similarly one may show that certain ramified
coverings of surfaces give rise to inclusions between mapping
class groups of closed surfaces. However, the intuition remains
that the closed surface mapping class groups are in some definite
sense as distinct, for differing values of $g$,  as possible.

A recent result of M.R. Bridson and K. Vogtmann on the analogous
question for automorphism groups of free groups shows that there
is no nontrivial morphism from the group Aut~$F_n$ into Aut~$F_m$
when $n>m>1$. In this paper we exploit some well-known facts about
Dehn twists, together with a specific twist description of a
torsion element of maximal order in the mapping class group
$\mcg_g$, to show that this same result holds true for $\mcg_g$.
In fact our methods show that only one or two genus-dependent
algebraic properties of these groups are needed, so that morphisms
into arbitrary groups lacking the relevant properties are
automatically trivial.

\section{Basic results on subgroups of mapping class groups}
We fix a closed oriented surface $S$ of genus $g$ as shown in
Figure~\ref{figure1}. Let $\mcg_g$ denote the mapping  class group
of $S$, the group of orientation preserving diffeomorphisms of
$S$ considered up to isotopy.

The Torelli group ${\mathcal{T}}_g$ is, by definition, the kernel
of the action of the mapping class group on the first homology
group $H_1(S;\Z)$ of $S$. The choice of a canonical basis for $H_1(S;\Z)$
gives rise to an epimorphism from the mapping class group $\mcg_g$
onto the symplectic group $Sp(2g,\Z)$ and one obtains a short
exact sequence
$$1\to {\mathcal{T}}_g \to \mcg_g \stackrel{\eta}\to Sp(2g,\Z) \to
1.$$

Notationally, the same letter will be used for a diffeomorphism and
its
isotopy class. Similarly, we make no distinction between isotopic
curves in a surface.

We introduce some basic types of mapping classes and state some
fundamental facts about them. For a simple closed curve $a$ on
$S$, we denote by $t_a$ the (isotopy class of a) right Dehn twist
along $a$. Two obvious simple properties of these mappings have
proved to be crucial in building up more intricate surface mappings
and also in recognising relations between them.
Firstly, Dehn twists along disjoint loops commute (up to isotopy):
this makes it
easy to find free abelian subgroups of rank $3g-3$ in $\mcg_g$, by
choosing maximal systems of disjoint homotopically distinct loops
which determine a pair of pants decomposition of $S$. Secondly, if
$f$ is any surface homeomorphism, then (again, up to isotopy)
$t_{f(a)}=f\circ t_a\circ
f^{-1}$.

It is elementary that if $a$ is a separating (i.e. null-homologous)
loop then
the Dehn twist $t_a$ is contained in the Torelli group
${\mathcal{T}}_g$.
Furthermore, if $a$ and $b$ are two non-separating simple closed
curves whose union bounds a subsurface of $S$,  then $t_at_b^{-1}$
is also in the group ${\mathcal{T}}_g$;  moreover, it follows from
Theorem~$2$ in~\cite{j} that if one of the two subsurfaces
bounded by the union $a\cup b$ is of genus $1$,
then ${\mathcal{T}}_g$ is the normal closure of
$t_at_b^{-1}$ in $\mcg_g$.  A precise
converse result by Vautaw will be discussed and used in the last
section.

We shall also need to focus attention on certain torsion elements
of $\mcg_g$. According to the classical Hurwitz-Nielsen Realisation
Problem, any finite subgroup of $\mcg_g$ should be isomorphic to a group
of automorphisms of some compact Riemann surface of genus $g$, and
a theorem of Hurwitz states that the order of such a group is at
most $84(g-1)$. The  Realisation Problem was proved by
S.P. Kerckhoff~\cite{ke}, but the case of finite cyclic groups had
been settled much earlier by J. Nielsen and W. Fenchel. This
characterisation of the finite cyclic subgroups of the mapping
class group by surface mappings leads to a bound on their order.
We state the basic facts we shall need about torsion subgroups
below ~\cite{w,h,h2}. For more details on automorphisms of Riemann
surfaces, the reader might consult  \cite{b,fk}.

\begin{theorem} \label{thmtorsion}
\item$(a)$  The order of any finite subgroup of $\mcg_g$ is at most
$84(g-1)$ if $g\geq 2$.
\item$(b)$ The order of a finite cyclic subgroup of
$\mcg_g$  is at most $4g+2$. This bound is achieved in every
genus.
\item$(c)$ If $\mcg_g$ has an element of prime order $p$,
then either $p\leq g+1$ or $p=2g+1$.
\item$(d)$ There is
no element of order $4g+1$ in $\mcg_g$.
\end{theorem}

Certain geometrically defined torsion elements, including a
specific element of maximal
order in $\mcg_g$, first documented by Wiman \cite{w} in 1895, will
feature prominently in the sequel.

\section{Homomorphisms between mapping class groups}

As background reference and source for unexplained facts,
we refer the reader to Ivanov's survey paper
\cite{iv}.

\begin{figure}[hbt]
 \begin{center}
   \includegraphics[width=10cm]{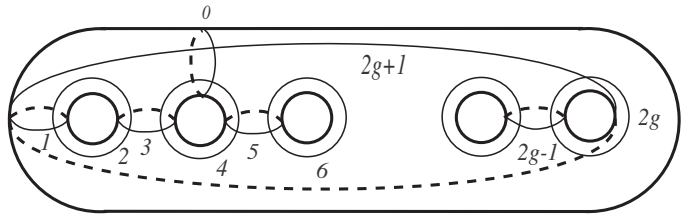}
  \caption{A system of simple loops in $S$.\hfill\break
  The curve labelled $i$ denotes $a_i$.}
  \label{figure1}
   \end{center}
 \end{figure}

Consider the system of simple closed curves
(loops) $a_0, a_1,a_2,\ldots, a_{2g+1}$ on $S$ as shown in
Figure~\ref{figure1} and write $t_j$ for the twist along the curve
$a_j$.  We let $\delta=t_1 t_2\cdots t_{2g}$. Note that for the
composition of two functions $\alpha$ and $\beta$, we use the
following notation: $\alpha\beta$ means $\alpha\circ \beta$, so
that $\beta$ is applied first.

\begin{lemma} \label{lemma1}
The order of $\delta$ in $\mcg_g$ is $4g+2$.
\end{lemma}
\begin{proof} Using the elementary properties of twist maps
mentioned
above, it is easy to verify that $\delta (a_i)=a_{i+1}$ for $1\leq
i \leq 2g$ and $\delta (a_{2g+1})=a_1^{-1}$: recall that equality
for closed loops in $S$  means they are freely
homotopic. It follows that the
order of $\delta$ is at least $4g+2$. On the other hand, from
Theorem~\ref{thmtorsion} we know that the order is at most $4g+2$.
The fact that $\delta^{4g+2}=1$ can also be proved by a more
elaborate direct computation with twists.
\end{proof}

\begin{lemma} \label{lemma1.5}
Let  $g\geq 1$ and let $a$ and $b$ be two non-separating simple
closed curves on $S$. Let $N$ be the normal closure of
$t_a^{-1}t_b$ in $\mcg_g$. Suppose either that $a$ intersects $b$
transversely at one point, or that $a$ is disjoint from $b$ and the
complement of $a\cup b$ in $S$ is connected. Then $N$ is
the commutator subgroup $[\mcg_g, \mcg_g]$. In particular,
$N=\mcg_g$ for $g\geq 3$.
\end{lemma}
\begin{proof}
The proof follows similar lines to a
result of McCarthy and Papadopoulos on involutions in $\mcg_g$,
\cite{mp}.  It is well-known that for $g\geq 2$ the mapping class
group
$\mcg_g$ is generated by the Dehn twists $t_0,t_1,\ldots, t_{2g}$
about the curves $a_0,a_1,\ldots,a_{2g}$ of Figure 1. When $g=1$,
the curve
$a_0$ does not exist, and $\mcg_1$ is generated by the two twists
$t_1, t_2$. We denote by $\Gamma$ the quotient group $\mcg_g/N$.

Suppose first that $a$ intersects $b$ at one point.  Clearly, for each
$i=1,2,\ldots, 2g-1$  there is a
diffeomorphism $f_i:S\to S$ such that $f_i(a)=a_i$ and $f_i(b)=a_{i+1}$.
Hence, since $t_a^{-1}t_b$ is in $N$,
so is the element
$$t_i^{-1}t_{i+1}=t_{f_i(a)}^{-1}t_{f_i(b)}=f_it_{a}^{-1}t_{b}f_i^{-1}.$$
Similarly, one sees that $t_0^{-1}t_4$
belongs to $N$. Consequently, all generators of $\mcg_g$
represent the same element in $\Gamma$, which implies that the
group $\Gamma$ is cyclic, hence abelian. Therefore, $N$ contains
$[\mcg_g, \mcg_g]$. On the other hand, since both $a$ and $b$
are non-separating, there exists a diffeomorphism $k:S\to S$
such that $k(a)=b$. Then  $t_a^{-1}t_b=t_a^{-1}k t_a k^{-1}$ is a
commutator,
which shows that $N$ is contained in $[\mcg_g, \mcg_g]$.
Thus we have  $N=[\mcg_g, \mcg_g]$.

Suppose next that $a$ is disjoint from $b$ and that the complement
of
$a\cup b$ is connected. In this case, $g$ must be at least $2$, of
course.
By cutting $S$ open along $b$, it is easy to see that there exists
a non-separating
simple loop $c$ on $S$ such that $c$ intersects $a$
transversely at one point, $c$ is disjoint from $b$ and the
complement of $b\cup c$ is connected. By the
classification of surfaces, there is a diffeomorphism
$f:S\to S$ such that $f(a)=b$ and $f(b)=c$. From this, it follows
that
$$t_a^{-1}t_c=t_a^{-1}t_b\,
t_b^{-1}t_c=t_a^{-1}t_bft_a^{-1}t_bf^{-1}$$
is contained in $N$. Now the conclusion of the lemma follows from
the first case.

For the final statement, we have $g\geq 3$, so that the first
homology group $H_1(\mcg_g;\Z)$
of $\mcg_g$ is trivial by Powell's theorem. This means that
$[\mcg_g , \mcg_g]$, which coincides with $N$, is equal
to $\mcg_g$.
\end{proof}

Our first theorem determines the normal closure of each
power of $\delta$ in $\mcg_g$.

\begin{theorem} \label{thmmain}
$(a)$ Let $k$ be an integer with $1\leq k\leq 2g$. The normal
subgroup of $\mcg_g$ generated by $\delta^k$ is the full group
$\mcg_g$ if $g\geq 3$. It has index $2$ if $g=2$ or $(g,k)=(1,1)$
and index $4$ if $(g,k)=(1,2)$.

$(b)$ If $g\geq 3$, then the normal subgroup of $\mcg_g$
generated by $\delta^{2g+1}$ contains the Torelli group
$\mathcal{T}_g$ as a subgroup of index $2$. In fact it is the
kernel of the natural morphism $\psi:\mcg_g\to PSp(2g,\Z)$, given
by
$\psi= P\circ \eta$ where $P: Sp \to PSp=Sp/\{\pm I\}$.
\end{theorem}

\begin{proof}
For each positive integer $k\leq 2g$, let $N_k$ denote the normal
subgroup of $\mcg_g$ generated by $\delta^k$. We prove first that
each $N_k$ contains the commutator subgroup of $\mcg_g$.

We observe that for any $k$ the element
$$t_1^{-1}\delta^k
t_1\delta^{-k}=t_1^{-1}t_{\delta^k(a_1)} =t_1^{-1}t_{k+1}$$ is
contained in $N_k$. Next we note that $a_1$ intersects
$a_{k+1}=\delta^k (a_1)$ at one point if $k=1$ or $k=2g$, whereas
if $2\leq k\leq 2g-1$ then $a_1$  is disjoint from $a_{k+1}$ and
the complement of $a_1\cup a_{k+1}$ in $S$ is connected. By
Lemma~\ref{lemma1.5}, $N_k$ contains the commutator subgroup of
$\mcg_g$, the normal closure of $t_1^{-1}t_{k+1}$. It follows that
the quotient $\mcg_g/N_k$ is abelian and hence a quotient of
$H_1(\mcg_g;\Z)$.

For $g\geq 3$, the
first homology group $H_1(\mcg_g;\Z)$ of $\mcg_g$ is trivial by
Powell's theorem. Therefore, $N_k=\mcg_g$ in this case.
If $g=2$, the group $H_1(\mcg_2;\Z)$ is isomorphic to $\Z/{10}\Z$
by a result of Mumford and is generated by the class of
the Dehn twist about any non-separating simple closed curve. The
element $\delta^k$ represents the $4k$-th power of this generator
of $H_1(\mcg_2;\Z)$.  Hence the quotient group $\mcg_2/N_k$ is
cyclic of order $2$. If $g=1$, then $H_1(\mcg_1;\Z)$ is isomorphic
to $\Z/{12}\Z$, again generated by the class of the Dehn twist
about any non-separating simple closed curve. The element
$\delta^k$ represents the $2k$-th power of this generator of
$H_1(\mcg_1;\Z)$, so it follows that $\mcg_1/N_k$ is a cyclic group,
of order $2$ if $k=1$ and of order $4$ if $k=2$. The conclusion
(a) of the theorem  follows.

To prove $(b)$, let $N$ denote the normal closure of
$\delta^{2g+1}$ in $\mcg_g$. It can easily be shown (e.g. by
another twist calculation) that $\delta^{2g+1}$ is  (up to
isotopy) a hyperelliptic involution of $S$, taking each loop
$a_j,\,1\leq j\leq 2g$ to its inverse. Thus, $\delta^{2g+1}$ acts
as minus the identity on the first homology of $S$. Therefore, $N$
is contained in the kernel of $\psi : \mcg_g\to PSp(2g,\Z)$. On
the other hand, the element
\[t_0^{-1} \delta^{2g+1} t_0
\delta^{-(2g+1)}=t_0^{-1}t_{\delta^{2g+1}(a_0)}\] is contained in
$N$. Now, $a_0$ and ${\delta^{2g+1}(a_0)}$ are disjoint
non-separating curves whose union bounds a subsurface of genus one,
and it follows from Theorem~$2$ in~\cite{j} that
the Torelli group is the normal closure of such an element.
Hence, $\mathcal{T}_g$ is contained in $N$. But the hyperelliptic
involution $\delta^{2g+1}$ is not
in $\mathcal{T}_g$ and the index of $\mathcal{T}_g$ in the kernel
of $\psi$ is $2$. Therefore $N$ is equal to the kernel of
$\psi$, proving $(b)$.

This completes the proof of the theorem.
\end{proof}

\begin{remark}
If $g=2$, then since $\delta^5$ is the hyperelliptic involution,
which is central in $\mcg_2$, the normal closure of $\delta^5$ is
cyclic of order
$2$. Thus, Theorem~\ref{thmmain}(b) is false in this case. For
$g=1$,
of course, the Torelli group is trivial.
\end{remark}

 \begin{figure}[hbt]
 \begin{center}
   \includegraphics[width=12cm]{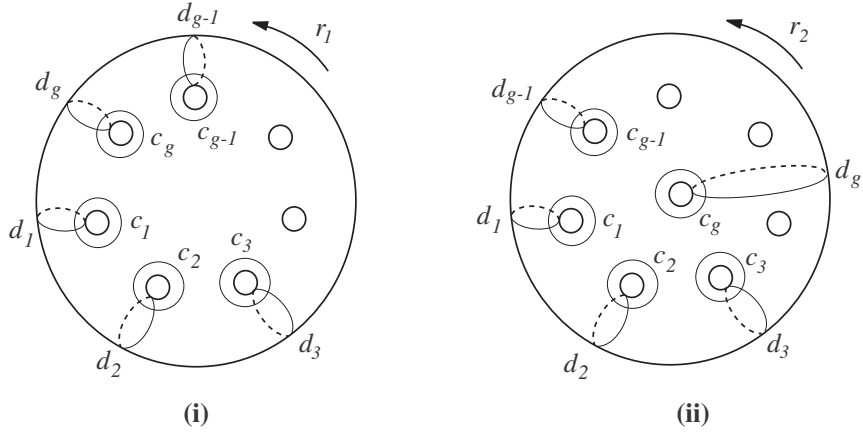}
  \caption{Two finite order surface mappings,
   induced by rotations $r_1$ and $r_2$ through $2\pi /g$ in
  $({\rm i})$ and $2\pi /(g-1)$ in $({\rm ii})$. }
  \label{figure2}
   \end{center}
 \end{figure}

\begin{lemma} \label{lemma-rot}
Let $g\geq 3$ and let $r$ denote either of the rotations $r_1$ and
$r_2$ illustrated in Figure~$\ref{figure2}$. Then the normal
closure of $r$ in $\mcg_g$ is the full mapping class group
$\mcg_g$.
\end{lemma}
\begin{proof}
In order to avoid double indices, we let $d$ denote the curve
$d_1$. Note that the curves $d$ and $r(d)=d_2$ are disjoint, the
complement of their union is connected and the element
$t_{d}^{-1}\, t_{r(d)}=t_{d}^{-1}\, r\, t_{d}\, r^{-1}$ is
contained in the normal closure of $r$. The lemma now follows from
Lemma~\ref{lemma1.5}, after applying a homeomorphism which identifies the surface
in Figure 1 with that of Fig 2; for instance we may take one which carries the curves $a_{2j-1}$
to $c_j$ for $j=1,\ldots,g$.
\end{proof}

We are now ready to prove our main result.

\begin{theorem} \label{triviality}
Let $g> h$ and let $\varphi : \mcg_g\to \mcg_h$ be a homomorphism.
The image of  $\varphi$ is trivial if $g\neq 2$ and has order at
most two if $g=2$.
\end{theorem}

\begin{proof}
There is nothing to prove when $h=0$, since the group $\mcg_0$ is
trivial.  Hence we assume that $h\geq 1$, so that $g\geq 2$.

The subgroup $\langle \delta \rangle$ of $\mcg_g$ generated by
$\delta$ is cyclic of
order $4g+2$ as we saw earlier. But by Theorem~\ref{thmtorsion}
(b), the maximal cyclic
subgroups of $\mcg_h$ have order at most $4h+2<4g+2$, and so there
must
be an integer $k$ with $1\leq k\leq 2g+1$ such that $\varphi
(\delta^k)=1$.  We remark also that if $g\geq 3$, the kernel of
$\varphi$ contains the Torelli group. Therefore, $\varphi$ factors
through the natural map $\eta:\mcg_g\to Sp(2g,\Z)$. That is,
$\varphi$ induces a map $\Phi : Sp(2g,\Z)\to \mcg_h$ such that
$\Phi\eta=\varphi$.

Now if the integer $k$, for which $\varphi (\delta^k)=1$, is less
than $2g+1$, it follows from Theorem~\ref{thmmain} that the
image of $\varphi$ is trivial if $g\geq 3$ and has order at most
$2$ if $g=2$, because the normal closure of $\delta^k$ is
contained in the kernel of $\varphi$.

The theorem follows immediately from this observation in the
(simplest possible) case
where $2g+1$ is prime. For since $2g+1$ is greater than $2h+1$, by
Theorem~\ref{thmtorsion}~$(c)$ the mapping class group $\mcg_h$ has
no element of order $2g+1$. Therefore we must have
$\varphi (\delta^k)=1$ for some $1\leq k\leq 2g$. This completes
the proof for infinitely many values of $g$, in particular for
$g=2,3$ and $5$.

Suppose next that $g> 5$. Consider a symplectic basis
$c_1,d_1,\ldots,c_g,d_g$ of $H_1(S;\Z)$ (with appropriate choice
of orientations) as shown in Figure~\ref{figure2}~(i) if $g$ is
odd and Figure~\ref{figure2}~(ii) if $g$ is even. There is a
subgroup $\Sigma_g$ of $Sp(2g,\Z)$, isomorphic to the symmetric
group on $g$ letters, which represents all permutations of the
hyperbolic pairs $(c_i,d_i)$. Let $A_g$ denote the alternating
subgroup of $\Sigma_g$. Let $r$ denote $r_1$ if $g$ is odd and
$r_2$ if $g$ is even in Figure~\ref{figure2}. Then the image
$\eta(r)$ of the surface mapping $r$ is contained in $A_g$.

The group $\mcg_h$ cannot contain a subgroup of order $g!/2$ when
$g>5$,  because $g!/2$ is greater than $84(h-1)$. Therefore, the
restriction of $\Phi$ to $A_g$ is not injective, and since $A_n$
is simple for $n> 5$, $A_g$ must be contained in the kernel of
$\Phi$. Therefore $r$ is contained in the kernel of $\varphi$ and,
since the normal closure of $r$ is the full group $\mcg_g$,
$\varphi$ is trivial.

Suppose finally that $g=4$. If $h=1$ then $4h+2$ is less than
$2g+1$. If $h=2$ then $4h+1=2g+1$. Therefore, by
Theorem~\ref{thmtorsion} we have $\varphi(\delta^k)=1$ for some
$1\leq k\leq 2g$. The conclusion of the theorem follows in this
case. Now let $h=3$. Again, $\varphi$ factors through the map
$\eta$. Consider the element $\gamma=t_1t_2t_3t_4$. A twist
calculation shows that $\gamma^{10}$ is the Dehn twist about a
nullhomologous simple closed curve, the boundary of $S$ cut open
along all these curves. Hence it is contained in the Torelli
group. Moreover, the order of $\eta(\gamma)$ is $10$ in
$Sp(8,\Z)$. Since the prime $5$ is greater than $h+1$ and
different from $2h+1$, there must be a positive integer $k<5$ such
that $\eta(\gamma)^k$ is contained in the kernel of $\Phi$, and so
$\gamma^k$ is contained in the kernel of $\varphi$. By considering
the curves $a_1$ and $\gamma^k(a_1)=a_{k+1}$, as in the proof of
Lemma~\ref{lemma1.5}, we see that the normal closure of $\gamma^k$
is the full group $\mcg_4$. Therefore, the image of $\varphi$ is
trivial in this case too.

This completes the proof.
\end{proof}

\begin{remark}
The method of proof for the general case $g>5$ in the above
theorem is a sharpening of our original one and was suggested to
us by a conversation with M.R. Bridson.
\end{remark}

\begin{corollary}
Let $g>5$ be an integer. Suppose that $H$ is any group with the
property that it contains no element of order $4g+2$ and no
subgroup of order $g!$.  Then any homomorphism $\varphi :
\mcg_g\to H$ has trivial image.
\end{corollary}
\begin{proof}
The proof of Theorem~\ref{triviality} applies verbatim.
\end{proof}

\section{Homomorphisms from mapping class groups}

In this section we complete the paper by drawing further
implications about arbitrary homomorphic images of mapping class
groups. The first amounts to no more than an abstraction of
Theorem~\ref{triviality}.

\begin{theorem} \label{2g+1}
Let $g$ be a positive integer. Suppose that $H$ is a group with
the property that it has no element of order $2g+1$.
If $\varphi : \mcg_g\to H$ is a homomorphism, then the image of
$\varphi$ is trivial if $g\geq 3$ and is a subgroup of a cyclic
group of order $2$ (resp. $4$) if $g=2$ (resp. $g=1$).
\end{theorem}

\begin{proof} The hypotheses imply that
there is an integer $k$ with $1\leq k\leq 2g$ such that $\varphi
(\delta^k)=1$. The normal closure $N_k$ of $\delta^k$ is contained
in the kernel of $\varphi$. But by
Theorem~\ref{thmmain}, it is the full group $\mcg_g$ if $g\geq 3$,
whereas it has index $2$ when $g=2$ and index $2$ or $4$ if $g=1$.
Since the quotient $\mcg_g/N_k$ is cyclic, the result follows.
\end{proof}

We note that if $g=2$ (resp. $g=1$), then there {\em is} a
homomorphism from the mapping class group $\mcg_g$ onto a cyclic
group of order $2$ (resp. $4$).

By the same method, but bringing into play the rotations $r_i$, we
obtain the next result.

\begin{theorem} \label{r_i}
Let $g\geq 3$ be a positive integer.
If $\varphi : \mcg_g\to H$ is a homomorphism, then the image of
$\varphi$ is trivial unless $H$ contains elements of order $g-1$
and $g$.
\end{theorem}
\begin{proof}
Suppose that $H$ contains no element of order $g$.
Consider the rotation $r_1$ in Figure~\ref{figure2}~(i).
Since the order of $r_1$ is $g$, there must be an integer $k$
with $1\leq k<g$ such that $\varphi(r_1^k)=1$. The proof of
Lemma~\ref{lemma-rot} also shows that the normal closure
of $r_1^k$ is equal to $\mcg_g$.
It follows that the kernel of $\varphi$ is $\mcg_g$.

If $H$ has no element of order $g-1$, the same argument applies,
mutatis mutandis, with the rotation $r_2$ of
Figure~\ref{figure2}~(ii) in place of $r_1$.
\end{proof}

A further criterion for recognising homomorphic images of $\mcg_g$
stems from restrictions on the rank of free abelian subgroups.
Recall that if $a$ is null-homologous then
$t_a\in{\mathcal{T}}_g$, and that if $b$ and $c$ are two disjoint
non-separating simple closed curves, whose union bounds a
subsurface of $S$,  then $t_b^{-1}t_c$ is in the Torelli group
${\mathcal{T}}_g$. A converse criterion, concerning the
intersection of any free abelian twist subgroup with the Torelli
group, is proved by Vautaw in~\cite{v}, Theorem~$3.1$. His result
can be stated as follows.

\begin{theorem} \label{thm:bill}
Let $c_1,c_2,\ldots,c_{k}$ be a set of homotopically nontrivial,
pairwise disjoint, non-separating,
simple closed curves on $S$ such that for each pair
$i\neq j$ the surface obtained by cutting
$S$ along $c_i$ and $c_j$ is connected.
Then an element $t_{c_1}^{n_1}t_{c_2}^{n_2}\cdots t_{c_k}^{n_k}$
is contained in the Torelli group
if and only if $n_i=0$ for all $i$.
\end{theorem}

As mentioned in section 1, it follows immediately from the definitions that Dehn twists
about homotopically  distinct loops commute up to isotopy, providing
a natural source of free abelian subgroups, though not the
only one.  There is, however, an upper bound for the rank of an
abelian subgroup of the mapping class group, determined by Birman,
Lubotzky and McCarthy~\cite{blm}.

\begin{theorem} \label{thm:blm}
The rank of a free abelian subgroup of the mapping class group
$\mcg_g$ is at most $3g-3$.
\end{theorem}

Certain specific subgroups of finite index in $\mcg_g$ will also
be significant here. We recall that, for any positive
integer $m$, the {\em congruence subgroup} $N(2g,m)$ of level $m$
in the genus $g$ symplectic group $Sp(2g,\Z)$ is the kernel of the
natural epimorphism from $Sp(2g,\Z)\to Sp(2g,\Z/m\Z)$ given by the
mod~$m$ reduction. The congruence subgroup problem for $Sp(2g,\Z)$
was solved by Mennicke~\cite{m}.

\begin{theorem} \label{thm:sp}
If $G$ is a normal subgroup of $Sp(2g,\Z)$ distinct from the
trivial subgroup and the center, then $G$ contains a congruence
subgroup $N(2g,m)$ for some $m$. In particular, $G$ is of finite
index in  $Sp(2g,\Z)$.
\end{theorem}
In the present context it is meaningful to ask whether the mapping
class groups possess a similar co-finality property in relation to
some natural infinite family of  finite index subgroups. However,
the range of possible subgroup families is considerably broader
than in the case of linear groups such as the symplectic modular
group and so far, in the absence of any faithful linear
representation for $\mcg_g$, there is no natural choice.

With this preamble we can state our final result about morphisms
from mapping class groups, which says essentially that in order
for a group to be an (infinite) homomorphic image of $\mcg_g$, it
must either have a torsion element of maximal order or the abelian rank
must be at least $3g-3$, that of $\mcg_g$.

\begin{theorem}
Let $g\geq 3$ be an integer. Suppose that $H$ is any group with the
following
properties: it contains no element of order $4g+2$ and any
free abelian subgroup has rank $<3g-3$. If $\varphi : \mcg_g\to H$
is a homomorphism, then the image of  $\varphi$ is finite.
\end{theorem}

\begin{proof}
Since $H$ has no element of order $4g+2$, there must be an integer
$k$ with $1\leq k\leq 2g+1$ such that $\varphi (\delta^k)=1$. If
$\varphi (\delta^k)=1$ for some $1\leq k\leq 2g$, then the image
of $\varphi$ is trivial by Theorem~\ref{thmmain} and the
conclusion is obvious in this case. If $\varphi (\delta^k)\neq 1$
for all $k\leq 2g$, then $\varphi (\delta^{2g+1})= 1$. But then
the kernel of $\varphi$ contains the Torelli group, the kernel of
$\eta$. Therefore there is a homomorphism $\Phi :  Sp(2g,\Z) \to
H$ such that $\Phi \eta = \varphi$.

It is an elementary exercise
to choose a system of non-separating, pairwise disjoint,
non-isotopic, simple closed curves $c_1,c_2,\ldots, c_{3g-3}$ on
$S_g$ such that
\item(i) the complement $S\setminus\cup_{j=1}^{3g-3} c_j$ is a
union of $2g-2$
spheres-with-three-holes, and
\item{(ii)} for each $i\neq j$ the complement of $c_i\cup c_j$ is
connected.

\noindent Now the (abelian) subgroup $A$ of $\mcg_g$ generated by the Dehn
twists about $c_1,c_2,\ldots, c_{3g-3}$ has rank $3g-3$. Also,
since the intersection $A\cap {\mathcal{T}}_g$ is trivial by
Theorem~\ref{thm:bill}, the restriction of $\eta$ to $A$ is
injective. Hence, $\eta (A)$ is a free abelian group of rank
$3g-3$. Since $\Phi(\eta(A))$ is abelian in $H$ and $H$ does not
contain any free abelian subgroup of
rank $3g-3$, there must be a
non-trivial element $h\in \eta(A)$ such that $\Phi(h)=1$. Clearly,
$\eta(A)$ does not contain the central element $-I$ in
$Sp(2g,\Z)$. Therefore, the kernel of $\Phi$ is a nontrivial
normal subgroup of $Sp(2g,\Z)$ different from the center and so,
by
Theorem~\ref{thm:sp}, it contains a congruence subgroup. In
particular, the image of $\Phi$ (and, hence, of $\varphi$) is
finite.

This completes the proof.
\end{proof}

\section{Concluding remarks}
Our results imply restrictions on the index of any
normal subgroup of finite index
in a mapping class group for $g\geq 3$: let  $G < \mcg_g$ be a
normal subgroup of finite index
$n$, then $n$ is divisible by $g-1$, $g$ and $2g+1$.
Results of this type may have some bearing on the notion
of subgroup growth for mapping class groups.

As N.V. Ivanov has pointed out, it is natural to ask a
more general question of rigidity type, whether a non-trivial
morphism (or more generally a morphism with infinite image)
exists from some {\em finite index subgroup of}
$\mcg_g$ to $\mcg_h$ with $g>h$. The abelian rank aspect of our
final result may have some bearing on answering this problem.
We note that there are subgroups of finite index in $\mcg_2$
admitting a homomorphism onto a free group of finite rank~\cite{ko}:
since $\mcg_1$ is virtually free, we must thus assume
in this problem that $g$ is at least $3$. It should also be made clear
that a negative answer to this problem would imply a positive answer
to Problem~$2.11$~$(A)$ in~\cite{ki}.

\bigskip
\noindent
{\bf Acknowledgements:} This paper grew from a discussion between
the first author and M.R. Bridson on the comparison between
Out\,$F_n$ and $\mcg_g$, prompting a question which the first author
posed in communication with Nikolai V. Ivanov, who included it in
a problem session at the AMS meeting in Ann Arbor in March 2002,
which the second author attended. Our joint work started when the second
author (M.K.) visited King's College London on 2002;  M.K  thanks
the Mathematics Department of King's College London for
hospitality and financial support. Both authors are grateful to
the Warwick Mathematics Institute for hospitality and to the
London Mathematical Society for travel to (and financial support
during) the Symposium on Geometric Topology at Warwick in July
2002.


\begin{thebibliography}{xxxx}

\bibitem[BLM]{blm}
J.~S.~Birman, A.~Lubotzky, J.~McCarthy,
{\em Abelian and solvable subgroups of the mapping class groups},
Duke Math. J. {\bf 50} (1983), 1107--1120.

\bibitem[B]{b} T.~Breuer,
{\em Characters and automorphism groups of compact Riemann
surfaces},
London Mathematical Society Lecture Note Series, 280.
Cambridge University Press, Cambridge, 2000.

\bibitem[BV]{bv}
M.~R.~Bridson, K.~Vogtmann, {\em  Homomorphisms from automorphism
groups of free groups}, To appear in Bull. London Math. Soc.

\bibitem[FM]{fm} B.~Farb, H.~Masur,
{\em Superrigidity of mapping class groups},
Topology {\bf 37} (1998), 1169-1176.

\bibitem[FK]{fk} H.~M.~Farkas, I.~Kra,
{\em Riemann Surfaces},
2nd ed., Graduate Texts in Mathematics, vol 71, Springer-Verlag,
Berlin Heidelberg New York, 1992.

\bibitem[H]{h} W.~J.~Harvey,
{\em Cyclic groups of automorphisms of a compact Riemann surface},
Quart. J. Math. Oxford Ser. (2) {\bf 17} (1966), 86--97.

\bibitem[H2]{h2} W.~J.~Harvey,
{\em On branch loci in Teichm\"uller space},
Trans. Amer. Math. Soc. {\bf 153} (1971), 387 - 399.

\bibitem[Iv]{iv} N.~V.~Ivanov,
{\em Mapping Class Groups}, Chapter 12 in Handbook  of Geometric
Topology, (Editors R.J. Daverman \& R.B. Sher), Elsevier Science
(2002), 523-633.

\bibitem[J]{j}
 D.~Johnson,
 {\em Homeomorphisms of a surface which act trivially on homology},
 Proc. Amer. Math. Soc. {\bf75} (1979), 119-125.

\bibitem[KM]{km}
V.~A. Kaimanovich, H.~Masur,
{\em The Poisson boundary of the mapping class group},
Invent. Math. {\bf 125} (1996), 221-264.

\bibitem[Ke]{ke}
S.~P.~Kerckhoff,
{\em The Nielsen realization problem},
Ann. of Math. (2) {\bf 117} (1983), 235--265.

 \bibitem[Ki]{ki}
 R.~Kirby, {\em Problems in low-dimensional topology},
 in Geometric Topology (W.~Kazez ed.) AMS/IP Stud. Adv. Math. vol 2.2,
 American Math.~Society, Providence 1997.

\bibitem[Ko]{ko}
M.~Korkmaz,
{\em On cofinite subgroups of mapping class groups},
Proceedings of 9th G\'okova Geometry-Topology Conference, and
Turkish Journal of Mathematics, to appear.

\bibitem[MP]{mp}
J.~D.~McCarthy, A.~Papadopoulos,
{\em Involutions in surface mapping class groups},
L'Enseignement Math\'ematique {\bf33} (1987), 275-290.

\bibitem[M]{m} J.~Mennicke,
 {\em Zur Theorie der Siegelschen Modulgruppe},
 Math. Ann. {\bf159} (1965), 115--129.

\bibitem[V]{v} W.~R.~Vautaw,
{\em Abelian subgroups of the Torelli group},
Alg. Geom. Topol. {\bf2} (2002), 157-170.

\bibitem[W]{w} A.~Wiman,
{\em Ueber die hyperelliptischen Curven und diejenigen vom
Geschlechte $p=3$,
welche eindeutigen Transformationen in sich zulassen},
Bihang Kongl. Svenska Vetenskaps-Akademiens Handl. (Stockholm
1895-6) Vol. {\bf 21}, 1-23.

\end{thebibliography}
\end{document}